# ASYMPTOTICALLY MINIMAX BAYES PREDICTIVE DENSITIES

By Mihaela Aslan

*Yale University*

Given a random sample from a distribution with density function that depends on an unknown parameter $\theta$, we are interested in accurately estimating the true parametric density function at a future observation from the same distribution. The asymptotic risk of Bayes predictive density estimates with Kullback–Leibler loss function $D(f_\theta||\hat{f}) = \int f_\theta \log(f_\theta/\hat{f})$ is used to examine various ways of choosing prior distributions; the principal type of choice studied is minimax. We seek asymptotically least favorable predictive densities for which the corresponding asymptotic risk is minimax. A result resembling Stein's paradox for estimating normal means by the maximum likelihood holds for the uniform prior in the multivariate location family case: when the dimensionality of the model is at least three, the Jeffreys prior is minimax, though inadmissible. The Jeffreys prior is both admissible and minimax for one- and two-dimensional location problems.

**1. Introduction.** There has been a historical dispute between the classical estimative density functions and the Bayesian predictive density functions in measuring the goodness-of-fit of the density estimate. For both the frequentist and Bayesian approaches to prediction inference, the choice of the prior is a serious matter either asymptotically or for finite samples. In this paper, we examine the asymptotic behavior of Bayes predictive density estimates under the Kullback–Leibler loss. These asymptotics are used to describe various ways of choosing prior distributions; the principal type of choice studied is minimax. Admissibility questions are also addressed for various families of densities.

Suppose we are given a random sample $\mathbf{x_n} = (x_1, x_2, \dots, x_n)$ of $n$ independent, identically distributed observations with respect to a probability density $f_\theta(\cdot) = f(\cdot|\boldsymbol{\theta})$, $\boldsymbol{\theta} \in \Theta \subseteq \mathcal{R}^p$, that depends on an unknown, $p$-dimensional

parameter. In the Bayesian approach, we assume some density $h(\cdot)$ over $\Theta$ to represent our prior knowledge of $\boldsymbol{\theta}$. A future observation $x_{n+1}$ from the same distribution is predicted by using a density $\hat{f}(\cdot|\mathbf{x_n})$, which is called a *predictive density*. We are interested in a density estimation problem where the actual parameter to be estimated is the density at the next observation $f(x_{n+1}|\boldsymbol{\theta})$, given the true, unknown parameter $\boldsymbol{\theta} = (\theta_1, \theta_2, \ldots, \theta_p)$.

A natural loss function used to measure the distance between the two densities for the next observation, $f_\theta$ and $\hat{f}$, is the Kullback–Leibler divergence,

$$D(f_\theta||\hat{f}) = \int \log \frac{f(x_{n+1}|\boldsymbol{\theta})}{\hat{f}(x_{n+1}|\mathbf{x_n})} f(x_{n+1}|\boldsymbol{\theta})\, dx_{n+1},$$

which is positive unless $f(x_{n+1}|\boldsymbol{\theta})$ coincides with $\hat{f}(x_{n+1}|\mathbf{x_n})$. This measure depends on $\boldsymbol{\theta}$ and the particular sample $\mathbf{x_n}$ observed. While not being a distance due to lack of symmetry, the Kullback–Leibler divergence produces standard results and consistent density estimates and, in general, leads to a more tractable problem than other loss functions (the $L^1$ distance for example). From the Bayesian point of view, the Kullback–Leibler loss has historically been the main tool for obtaining noninformative priors; Jeffreys [9] used its invariance properties and local behavior as a Euclidean square of a distance function as a starting point in constructing and proposing the prior that carries his name.

The Jeffreys prior density with respect to the $p$-dimensional Lebesgue measure,

$$J(\boldsymbol{\theta}) \propto \det{}^{1/2}((L_{ij}(\boldsymbol{\theta}))_{i,j=1,\ldots,p}),$$

where $(L_{ij}(\boldsymbol{\theta}))_{i,j=1,\ldots,p}$ is the information matrix $\mathbf{P}_{\boldsymbol{\theta}}[-\partial^2/\partial\theta_i\,\partial\theta_j \log f_{\boldsymbol{\theta}}]$ and $\mathbf{P}_{\boldsymbol{\theta}}$ represents the expectation with respect to $f_{\boldsymbol{\theta}}$, plays an important role in our framework. Inferences around the Jeffreys density are very suitable here, especially in invariance-related problems. It is asymptotically least favorable under entropy risk [5], and for $\alpha = \frac{1}{2}$, belongs to the family of relatively invariant priors proposed by Hartigan [7],

$$\left\{\frac{\partial \log h}{\partial\theta_j} = \sum_{i,r} L_{i,r}^{-1} \mathbf{P}_{\boldsymbol{\theta}}\left[\alpha \frac{\partial \log f_{\boldsymbol{\theta}}}{\partial\theta_j} \frac{\partial \log f_{\boldsymbol{\theta}}}{\partial\theta_i} \frac{\partial \log f_{\boldsymbol{\theta}}}{\partial\theta_r} + \frac{\partial^2 \log f_{\boldsymbol{\theta}}}{\partial\theta_j\,\partial\theta_i} \frac{\partial^2 \log f_{\boldsymbol{\theta}}}{\partial\theta_j\,\partial\theta_r}\right]\right\}_\alpha.$$

This family of prior densities is not equivalent to the family of all relatively invariant priors and we will refer to it as the $\alpha$-family (or $\alpha$-class).

The risk function is the expected Kullback–Leibler loss with respect to $f_{\boldsymbol{\theta}}$. We consider this to be our measure of the goodness-of-fit of $\hat{f}(x_{n+1}|\mathbf{x_n})$ to the unknown $f(x_{n+1}|\boldsymbol{\theta})$:

$$\begin{aligned} R(\boldsymbol{\theta}, \hat{f}) &= \mathbf{P}_{\boldsymbol{\theta}}(D(f_{\boldsymbol{\theta}}||\hat{f})) \\ &= \int\int \log \frac{f(x_{n+1}|\boldsymbol{\theta})}{\hat{f}(x_{n+1}|\mathbf{x_n})} f(x_{n+1}|\boldsymbol{\theta}) f(\mathbf{x_n}|\boldsymbol{\theta})\, dx_{n+1}\, d\mathbf{x_n}. \end{aligned}$$

We also consider as our density estimate the Bayes predictive density for the next observation based on the prior $h(\boldsymbol{\theta})$ and the data $\mathbf{x_n}$,

$$f_h(x_{n+1}|\mathbf{x_n}) = \int f(x_{n+1}|\boldsymbol{\theta}) h(\boldsymbol{\theta}|\mathbf{x_n})\, d\boldsymbol{\theta},$$

where $h(\boldsymbol{\theta}|\mathbf{x_n})$ is the posterior density obtained by using Bayes' product formula,

$$\frac{h(\boldsymbol{\theta}) f(\mathbf{x_n}|\boldsymbol{\theta})}{\int h(\boldsymbol{\theta}) f(\mathbf{x_n}|\boldsymbol{\theta})\, d\boldsymbol{\theta}}.$$

Thus, under the Kullback–Leibler loss, the predictive Bayes density estimate for the next observation is just the posterior density of the next observation.

For samples of finite size, Aitchison [1] shows that when a specific prior density $h(\boldsymbol{\theta})$ is given, any estimative density $\hat{f}(x_{n+1}|\mathbf{x_n})$ is inferior, in Kullback–Leibler risk, to the Bayes predictive density $f_h(x_{n+1}|\mathbf{x_n})$. From the asymptotic point of view, Komaki [10] gives an asymptotic expression for the Bayesian predictive distribution and shows that in the multidimensional curved exponential family case, the estimative distributions, given asymptotically efficient estimators, can be improved to predictive distributions that asymptotically coincide with the Bayesian predictive distributions.

Following the program of Hartigan [8] for finding the maximum likelihood prior density, we are searching for a prior distribution corresponding asymptotically to the minimax risk, as the number of observations $n$ from $f_\theta$ increases to $\infty$. We use asymptotic expansions of Kullback–Leibler risk functions in which the first order term $\frac{p}{2n}$ is the same for all estimative and Bayes predictive densities, given any continuously twice differentiable positive prior density; we allow prior densities to have infinite total mass $\int_\Theta h(\boldsymbol{\theta})\, d\boldsymbol{\theta} = \infty$.

Finding asymptotically minimax Bayes predictive density estimates is usually a hard task for general statistical settings, especially due to infinite parameter spaces. Choosing prior density functions for which the asymptotic risk is minimax mainly reduces to solving very complicated differential equations in many dimensions, which may or may not have solutions. Even in noninvariant settings, by concentrating on a smaller class of priors with useful invariance properties (such as a class of relatively invariant priors), these differential equations become much simpler and we are sometimes able to arrive at minimax solutions.

The main idea of this paper is to describe a searching algorithm for least favorable priors which starts by looking for minimax solutions among relatively invariant priors in the $\alpha$-class. We compare different predictive density estimates by looking at the smaller order terms in the asymptotic risk. These $\frac{1}{n^2}$ terms involve expressions in both the likelihood and the prior. Thus, choosing one density estimate over another reduces mainly to choosing prior density functions that improve on the asymptotic risk. Admissibility and minimaxity questions are expressed in terms of certain differential operators; the answers to these questions are then determined by the existence of solutions to different partial differential equations.

*Algorithm scheme.*

1. Compute asymptotic risk expressions of the form
$$\frac{A}{n} + \frac{B}{n^2} + \cdots:$$
   - $A$ is the same constant for all estimative and Bayes predictive densities.
   - Different density estimates compete through the $\frac{1}{n^2}$ term, which depends on the likelihood and the prior.
2. Find priors leading to asymptotically minimax density estimates:
   - Start by searching in a smaller class of priors and find those priors for which the asymptotic risk is constant in the parameter (in other words, find the optimum in a wide class of possible estimates).
   - Prove that the priors with the smallest constant risk are least favorable: show that they cannot be uniformly beaten over all priors by solving a differential equation in the parameter (in other words, show that this optimum is also the optimum among all possible estimates).

This method is not restricted to relatively invariant priors in the $\alpha$-class or to invariant statistical problems. The $\alpha$-class merely represents a good "set of guesses" for an optimal estimate in the minimax sense. The methodology can be generalized to various distribution functions and, hence, to general statistical settings which do not present any symmetries or invariance properties.

One important application of this method is to the general location model: a result resembling Stein's paradox for estimating normal means by the maximum likelihood holds for the uniform prior in the multivariate location family case. Using differential geometry and, in particular, potential theory, we show that when the dimensionality of the location model is at least three, Jeffreys' prior is minimax, though inadmissible. The Jeffreys prior is both admissible and minimax for one- and two-dimensional location problems.

## 2. General notation and main result.

We begin by introducing certain notation that will be used throughout our discussion. Let

$$L(\boldsymbol{\theta}) = \prod_{j=1}^{n} f_{\boldsymbol{\theta}}(x_j),$$

$$l(\boldsymbol{\theta}) = \log L(\boldsymbol{\theta}) = \sum_{j=1}^{n} \log f_{\boldsymbol{\theta}}(x_j)$$

be the likelihood and the log-likelihood functions for the sample of observations $\mathbf{x_n}$. Also, let

$$l_{\mathbf{i}} = l_{i_1,\ldots,i_r} = \frac{\partial^r}{\partial\theta_{i_1}\cdots\partial\theta_{i_r}} \log L(\boldsymbol{\theta}) = \sum_{j=1}^{n} \frac{\partial^r}{\partial\theta_{i_1}\cdots\partial\theta_{i_r}} \log f_{\boldsymbol{\theta}}(X_j),$$

$$L_{\mathbf{i}_1,\mathbf{i}_2,\ldots,\mathbf{i}_s} = \mathbf{P}_{\boldsymbol{\theta}}[l_{\mathbf{i}_1} l_{\mathbf{i}_2} \cdots l_{\mathbf{i}_s}]$$

be the log-likelihood derivatives and the expectations of their products, all evaluated at the true value $\boldsymbol{\theta}$; $\mathbf{P}_{\boldsymbol{\theta}}$ denotes expectation, given $\boldsymbol{\theta}$. The same quantities, evaluated at the maximum likelihood estimate $\hat{\boldsymbol{\theta}}$, will be denoted by $\hat{L}$, $\hat{l}$, $\hat{l}_{i_1,\ldots,i_r}$ and $\hat{L}_{\mathbf{i}_1,\mathbf{i}_2,\ldots,\mathbf{i}_s}$. The matrix $(-L_{ij})_{i,j=1,\ldots,p}$ is called the *Fisher information matrix.*

The model is supplemented by a prior distribution on $\Theta$ with density function $h$ with respect to the Lebesgue measure, where $h_i = \frac{\partial}{\partial\theta_i} \log h$ and $h_{ij} = \frac{\partial^2}{\partial\theta_i \partial\theta_j} \log h$ stand for the log prior first and second derivatives when evaluated at the true $\boldsymbol{\theta}$. Let $\hat{h}$, $\hat{h}_i$ and $\hat{h}_{ij}$ be the same quantities when evaluated at the maximum likelihood estimate.

The following theorem gives the asymptotic expression for the Kullback–Leibler risk up to smaller $\frac{1}{n^3}$ terms. Adopting tensor summation conventions, the various expressions that appear in our formula are in fact sums of terms over indices that appear twice.

THEOREM 1. *Under regularity conditions stated in the Appendix, the asymptotic risk with terms of order $O(n^{-1})$ and $O(n^{-2})$, and ignoring smaller terms of order $O(n^{-3})$, has the following expression:*

$$\begin{aligned}
R(\boldsymbol{\theta}, f_h) = {} & \frac{p}{2n} - \frac{p}{4n^2} \\
& + \frac{1}{n^2}\biggl[n\biggl\{L_{i,r}^{-1} L_{j,s}^{-1}\left(\frac{1}{2}L_{ij,r,s} + \frac{3}{4}L_{ij,rs} + L_{irj,s} + \frac{1}{2}L_{irjs}\right) \\
& \qquad + L_{i,r}^{-1} L_{j,s}^{-1} L_{k,t}^{-1}\biggl(\frac{1}{2}L_{i,rj}L_{k,st} + \frac{1}{2}L_{i,jk}L_{t,rs} + \frac{1}{6}L_{ijk}L_{r,s,t} \\
& \qquad\qquad + L_{irj}L_{k,st} + \frac{3}{2}L_{ijk}L_{r,st} \\
& \qquad\qquad\qquad + \frac{1}{2}L_{irj}L_{skt} + \frac{7}{12}L_{ijk}L_{rst}\biggr) \\
& \qquad\qquad + L_{i,r}^{-1} L_{j,s}^{-1}(L_{rj,s} + L_{rjs})h_i \\
& \qquad\qquad\qquad + L_{i,r}^{-1}\left(h_{ir} + \frac{1}{2}h_i h_r\right)\biggr\}\biggr] \\
& + O(n^{-3}).
\end{aligned}$$

*Here $p$ is the dimensionality of $\Theta$.*

REMARK 1. $L^{-1}_{i,r}$ denotes the $(i,r)$ element of the inverse of the Fisher information matrix $(-L_{ir})_{i,r=1,\dots,p}$.

REMARK 2. The expression $n\{\cdots\}$ is the same for all $n$.

REMARK 3. The first-order term in this expansion coincides with the same order term in the asymptotic risks of the maximum likelihood and Bayes procedures ([8], Theorems 1 and 4). It also coincides with the upper bound for the asymptotic entropy risk from [5].

REMARK 4. The prior expressions from the second order terms in the asymptotic risks of the Bayes predictive densities and of the Bayes estimators from [8] are the same. Thus, the difference in the $O(n^{-2})$ terms of the two asymptotic risks does not depend on the prior choice. The estimative density for maximum likelihood estimate is $f(x|\hat{\boldsymbol{\theta}})$. It may be shown that the ratio of predictive to estimative densities is asymptotically the same for all priors, therefore, the difference in the Kullback–Leibler distances does not depend on the prior and, indeed, the estimative density has risk no less than the predictive density.

REMARK 5. Let $J$ denote the Jeffreys prior density. By adding and subtracting terms involving log Jeffreys' first and second derivatives, $J_i = \frac{\partial}{\partial\theta_i}\log J$ and $J_{ij} = \frac{\partial^2}{\partial\theta_i\,\partial\theta_j}\log J$, all evaluated at the true $\boldsymbol{\theta}$, the $O(n^{-2})$ part of the asymptotic risk expression will separate into two distinct expressions, both invariant under monotone transformations of the parameter. We call these new expressions *the likelihood* and *the prior terms*:

$$
\begin{aligned}
&\text{likelihood term}\\
&\quad = L^{-1}_{i,r}L^{-1}_{j,s}\left(\tfrac{1}{2}L_{ij,r,s} + \tfrac{3}{4}L_{ij,rs} + L_{irj,s} + \tfrac{1}{2}L_{irjs}\right)\\
&\qquad + L^{-1}_{i,r}L^{-1}_{j,s}L^{-1}_{k,t}\left(\tfrac{1}{2}L_{i,rj}L_{k,st} + \tfrac{1}{2}L_{i,jk}L_{t,rs} + \tfrac{1}{6}L_{ijk}L_{r,s,t}\right.\\
&\qquad\qquad\qquad \left. + L_{irj}L_{k,st} + \tfrac{3}{2}L_{ijk}L_{r,st} + \tfrac{1}{2}L_{irj}L_{skt} + \tfrac{7}{12}L_{ijk}L_{rst}\right)\\
&\qquad + L^{-1}_{i,r}L^{-1}_{j,s}\left(\tfrac{1}{2}L_{ir,j,s} + 2L_{ij,r,s} + L_{ij,rs} + L_{irj,s} + \tfrac{1}{2}L_{i,r,j,s}\right)\\
&\qquad + L^{-1}_{i,r}L^{-1}_{j,s}L^{-1}_{k,t}\left(-2L_{rs,t}L_{i,j,k} - L_{ij,k}L_{rs,t} - L_{ij,k}L_{rt,s}\right.\\
&\qquad\qquad\qquad + \tfrac{1}{2}L_{ij,r}L_{st,k} - \tfrac{1}{2}L_{i,j,k}L_{r,s,t}\\
&\qquad\qquad\qquad\qquad \left. + \tfrac{1}{8}L_{i,j,s}L_{r,k,t} + \tfrac{1}{2}L_{rk,t}L_{i,j,s}\right)\\
&\qquad + L^{-1}_{i,r}L^{-1}_{j,s}(L_{j,rs} + L_{rjs})L^{-1}_{k,t}\left(\tfrac{1}{2}L_{i,k,t} + L_{ik,t}\right);
\end{aligned}
$$

prior term

$$= L^{-1}_{i,r}\{-J_r(h_i - J_i) + (h_{ir} - J_{ir}) + \tfrac{1}{2}(h_i - J_i)(h_r - J_r)\}.$$

REMARK 6. If the group of invariant transformations of the parameter is transitive (so that a transformation exists mapping any parameter value into any other), then the likelihood term in the risk is constant. In this case, for relatively invariant priors, the prior term is also constant. Thus, for priors in the class of relatively invariant priors proposed by Hartigan [7],

$$\{h_i = L^{-1}_{j,s}(\alpha L_{i,j,s} + L_{ij,s})\}_\alpha,$$

the asymptotic risk expression is independent of the parameter and reduces to a quadratic function in $\alpha$. Solving for $\alpha$, one finds that the choice of $\alpha$ giving the asymptotically minimum risk satisfies the following condition:

$$\begin{aligned}(\alpha - 1)&L^{-1}_{i,r}L^{-1}_{j,s}L^{-1}_{k,t}L_{i,j,s}L_{r,k,t}\\ &= L^{-1}_{i,r}L^{-1}_{j,s}(-L_{ir,j,s} - L_{i,jr,s} - L_{i,j,rs} - L_{i,r,j,s})\\ &\quad + L^{-1}_{i,r}L^{-1}_{j,s}L^{-1}_{k,t}(L_{r,js}L_{i,k,t} + 2L_{rs,k}L_{i,j,t} + L_{r,s,k}L_{i,j,t}).\end{aligned}$$

The Jeffreys prior corresponds to $\alpha = \frac{1}{2}$ and is a member of the class.

2.1. *The one-parameter problem with examples.* A simpler asymptotic risk expression holds when $\boldsymbol{\theta}$ is one-dimensional.

COROLLARY 1. *For $\theta \in \Theta \subseteq \mathcal{R}$, the asymptotic risk expression becomes*

$$\begin{aligned}R(\theta, f_h) = {}& \frac{1}{2n} - \frac{1}{4n^2}\\ &+ \frac{1}{n^2}\biggl[n\biggl\{L^{-2}_{1,1}\biggl(\frac{1}{2}L_{1,1,2} + \frac{3}{4}L_{2,2} + L_{1,3} + \frac{1}{2}L_4\biggr)\\ &\qquad + L^{-3}_{1,1}\biggl(L^2_{1,2} + \frac{1}{6}L_3L_{1,1,1} + \frac{5}{2}L_3L_{1,2} + \frac{13}{12}L^2_3\biggr)\\ &\qquad\quad + L^{-1}_{1,1}(L_{1,2} + L_3)h_1 + L^{-1}_{1,1}\biggl(h_2 + \frac{1}{2}h^2_1\biggr)\biggr\}\biggr]\\ &+ O(n^{-3}).\end{aligned}$$

Having precise expressions for the asymptotic risk permits detailed evaluations of admissibility and minimaxity. In the following subsections, we present some applications to both discrete and continuous distribution functions. Although it is true that some of these examples can be done in finite sample settings, our general technique agrees with them and offers a method of arriving at minimax solutions. To simplify our calculations, we will assume that $n = 1$.

2.1.1. *The Poisson example.* Let $x$ be an observation according to the Poisson distribution. From Corollary 1, the $O(\frac{1}{n^2})$ term in the asymptotic risk that depends on the parameter is of the form

$$h_1 + \theta\Big(h_2 + \frac{1}{2}h_1^2\Big) - \frac{1}{12\theta}.$$

Within the class of relatively invariant priors $\{h = \theta^{\alpha-1}\}_\alpha$, the priors $h = \theta^{\pm 1/\sqrt{6}}$ corresponding to $\alpha = 1 \pm \frac{1}{\sqrt{6}}$ have constant risk. The prior corresponding to $\alpha = 1$ has risk everywhere smaller than these priors, so they are inadmissible. However the maximum risk for any prior is never less than the risk for these priors, so they are minimax.

2.1.2. *The binomial example.* Let $x$ be an observation according to the binomial distribution with the canonical parameter $\theta$, $Bin(1, \frac{e^\theta}{1+e^\theta})$. As in the Poisson case, the $\frac{1}{n^2}$ term of the asymptotic risk depending on the parameter is easily computed as being

$$\frac{1}{e^\theta}\Big\{\Big(h_2 + \frac{1}{2}h_1^2\Big)(1+e^\theta)^2 - h_1(1-e^{2\theta}) + \frac{5}{24} + \frac{1}{12}e^\theta + \frac{5}{24}e^{2\theta}\Big\}.$$

It can be shown that the prior corresponding to $\alpha = 1 + \frac{1}{\sqrt{6}}$ has constant risk and is minimax among all positive priors.

2.1.3. *The negative binomial example.* For the negative binomial distribution in the canonical parameter $\theta$, $\mathcal{N}Bin(r, 1-e^\theta)$, the $\frac{1}{n^2}$ term in the asymptotic risk is of the form

$$\frac{1}{re^\theta}\Big\{\Big(h_2 + \frac{1}{2}h_1^2\Big)(1-e^\theta)^2 - h_1(1-e^{2\theta}) + \frac{5}{24} - \frac{1}{12}e^\theta + \frac{5}{24}e^{2\theta}\Big\}.$$

Following the binomial case, the prior corresponding to $\alpha = 1 - \frac{1}{\sqrt{6}}$ gives least constant risk within the $\alpha$-family; we have not been able to show minimaxity.

2.1.4. *The normal location-scale example.* Suppose we have an observation $x$ according to the normal location-scale density function, $\mathcal{N}(\mu, \sigma^2)$. Due to obvious invariances with respect to groups of transformations over the sample and the parameter spaces, the asymptotic risk expression reduces to

$$R(\sigma^2, f_h) = \frac{1}{2n} + \frac{1}{n^2}\Big\{-\frac{1}{2} + \frac{19}{12} + 4\sigma^2 h_2 + 2\sigma^4\Big(h_{22} + \frac{1}{2}h_2^2\Big) + \frac{3}{4}\Big\} + O(n^{-3}),$$

where each subscript 2 for the prior represents differentiation with respect to $\sigma^2$.

Unlike the normal scale case, the Jeffreys prior, $J(\mu, \sigma^2) \propto \sigma^{-3}$, is neither admissible nor minimax: the prior corresponding to $\alpha = \frac{2}{3}$ in the $\alpha$-class has a strictly smaller asymptotic risk than Jeffreys' for all $\boldsymbol{\theta}$. This agrees with the finite sample result that the Bayes predictive density based on $\sigma^{-1}\,d\mu\,d\sigma$ (which corresponds to our $\alpha = \frac{2}{3}$) is the best invariant predictive density and has a strictly smaller asymptotic risk than the Jeffreys prior [12].

2.1.5. *The multivariate normal scale example.* Let $\mathbf{x_n} = (x_1, \ldots, x_n)$ be a random sample according to the multivariate normal scale distribution with the log-likelihood function

$$l(V) = \log f(x|V) = -\tfrac{1}{2}\sum_{i,j=1}^{n} W_{ij}x_i x_j - \tfrac{1}{2}\log|V| + \mathbf{ct},$$

where $V$ is a symmetric and positive definite covariance matrix, $V = W^{-1}$. Also let

$$l_{(\mathbf{ij})} = l_{(i_1 j_1)\cdots(i_r j_r)} = \sum_{k=1}^{n} \frac{\partial^r}{\partial V_{i_1 j_1}\cdots\partial V_{i_r j_r}} \log f_V(X_k),$$

$$L_{(\mathbf{i}_1\mathbf{j}_1),(\mathbf{i}_2\mathbf{j}_2),\ldots,(\mathbf{i}_s\mathbf{j}_s)} = \mathbf{P}_V\,[l_{(\mathbf{i}_1\mathbf{j}_1)} l_{(\mathbf{i}_2\mathbf{j}_2)}\cdots l_{(\mathbf{i}_s\mathbf{j}_s)}],$$

be the log-likelihood derivatives and the expectations of products of the log-likelihood derivatives, where single indices represent pairs of indices identifying the variance–covariance parameters. For example, $\operatorname{var}(X_i) = V_{ii} = W_{ii}^{-1}$, $\operatorname{cov}(X_i, X_j) = V_{ij} = W_{ij}^{-1}$. For each pair of indices $(ii')$, it is assumed that $i \le i'$.

We give, without proof, the following lemma that can be found in [2]:

LEMMA 1. *For any pairs of indices $(ii')$, $(rr')$, $(jj')$, $(ss')$, the following expressions for the expectations of different products of the log-likelihood derivatives are true:*

$$L_{(ii')} = 0,$$

$$L_{(ii')(rr')} = -\frac{W_{ir}W_{i'r'} + W_{ir'}W_{i'r}}{2^{\{i=i'\}+\{r=r'\}}},$$

$$\begin{aligned} L_{(ii')(rr')(jj')} &= \frac{1}{2^{\{i=i'\}+\{r=r'\}+\{j=j'\}}} \\ &\quad\times\{W_{i'j}W_{r'i}W_{j'r} + W_{ij}W_{r'i'}W_{j'r} + W_{i'r}W_{j'i}W_{r'j} \\ &\qquad + W_{ir}W_{j'i'}W_{r'j} + W_{i'j}W_{ri}W_{j'r'} + W_{i'r'}W_{j'i}W_{rj} \\ &\qquad + W_{i'j'}W_{r'i}W_{jr} + W_{i'r}W_{ji}W_{r'j'}\}, \end{aligned}$$

$L_{(ii')(rr')(jj')(ss')}$

$$= -\frac{1}{2^{\{i=i'\}+\{r=r'\}+\{j=j'\}+\{s=s'\}}}$$

$$\times \{\textit{sum of } 48 \textit{ terms like } W_{is}W_{i'j}W_{j'r}W_{r's'},$$

*where, for each pair of $W$'s, the indices must come from at least three pairs of indices in the $L$'s. For example, $W_{is}W_{i'j}$ comes from the pairs $ii', jj', ss'$*$\}$.

*Also, the inverse information matrix components are*

$$L^{-1}_{(ii'),(rr')} = V_{ir}V_{i'r'} + V_{ir'}V_{i'r}.$$

Using Remark 6 above, in the $\alpha$-class of priors, which is of the form

$$\{h_{ii'} = \alpha L^{-1}_{(jj')(ss')} L_{(ii'),(jj'),(ss')}\}_\alpha,$$

the choice of $\alpha$ giving the asymptotically minimum risk satisfies the following condition:

$$\begin{aligned}(\alpha-1)&L^{-1}_{(ii'),(rr')}L^{-1}_{(jj'),(ss')}L^{-1}_{(kk'),(tt')}L_{(ii')(jj')(ss')}L_{(rr')(kk')(tt')}\\ &= L^{-1}_{(ii'),(rr')}L^{-1}_{(jj'),(ss')}(L_{(ii')(rr')(jj')(ss')}\\ &\qquad + L_{(ii'),(rr')}L_{(jj'),(ss')} - L_{(rr'),(ss')}L_{(ii'),(jj')})\\ &\quad + L^{-1}_{(ii'),(rr')}L^{-1}_{(jj'),(ss')}L^{-1}_{(kk'),(tt')}L_{(rr')(ss')(kk')}L_{(ii')(jj')(tt')}.\end{aligned} \tag{2.1}$$

Through simple manipulations of the likelihood identities (A5) in the Appendix and the formulae in Lemma 1, explicit expressions for all the terms in (2.1) become available. Due to invariance arguments, it can be shown that the $\alpha = \frac{1}{2}$ solution to (2.1), which corresponds to the Jeffreys prior, has the minimum asymptotic risk within the $\alpha$-class of priors.

In finite sample theory, Murray [13] and Ng [14] prove similar results for general group models under invariant prediction.

For the univariate normal scale case, the Jeffreys prior is also minimax among all smooth priors available: through a simple reparametrization of the form $u = \log \sigma^2$, and following the argument in Section 3, this problem becomes a location problem for which one can prove that the uniform prior in the new parameter is least favorable among all smooth priors, and so is minimax.

**3. Minimaxity and admissibility in the location case.** One important application of the searching method for minimax solutions is to the general location model. As in the Stein problem of estimating multivariate normal location parameters, we prove that a similar division between dimensions 2 and 3 holds for predictive density estimates for the general location problem. For dimensions 1 and 2, the Jeffreys prior is admissible, but for dimensions greater than 2 there are priors that have everywhere smaller risk than Jeffreys' so that the Jeffreys prior, though minimax, is inadmissible.

In finite sample theory, the Stein phenomenon in density estimation has been explored by Komaki [11] who showed that for the multivariate normal location model, the Jeffreys prior produces density estimates admissible in 1 or 2 dimensions, but inadmissible in 3 or more, just as Stein did for location estimates. For the same multivariate normal location model, George, Liang and Xu [6] go further than Komaki and show that under certain conditions on the marginal of the prior, the corresponding Bayes predictive density becomes minimax.

We prove a more general result by using the asymptotic risk expression from Theorem 1 to evaluate admissibility and minimaxity of density estimates in a general location model setting. The risk evaluations require the study of elliptic differential operators. For parameter estimation, such operators, in a simple form, appear in [4]. We give here the general form of such operators for density estimates.

Let $f(\mathbf{x}-\boldsymbol{\mu})$ represent a general location model with the standard probability density function $f(\mathbf{x})$ and $\boldsymbol{\mu}$ the location parameter for the family. For general multivariate location families, the likelihood term is constant under invariant and transitive transformations of the parameter. Thus, the risk expression depends on the parameter only through the prior term, which has the expression

$$\sum_{i,r=1}^{p} L_{i,r}^{-1}\{-J_r(h_i-J_i)+(h_{ir}-J_{ir})+\tfrac{1}{2}(h_i-J_i)(h_r-J_r)\},$$

where $h$ is any continuously twice differentiable positive prior among all priors available and $J$ is the Jeffreys prior. The Jeffreys prior being constant, all of the log-Jeffreys derivatives are 0. Note that for $h=J$, the prior term is 0.

Choosing $h$ of the form $g^2$, $g>0$, the part that remains to be optimized becomes

$$\sum_{i,r=1}^{p} L_{i,r}^{-1}\{g_{ir}+g_i g_r\}=\sum_{i,r=1}^{p} L_{i,r}^{-1}\frac{\partial^2/\partial\mu_i\,\partial\mu_r g}{g}=\frac{\Delta g}{g},$$

where $\Delta$ stands for the Laplacian differential operator and $\Delta g=\sum_{i=1}^{p}\frac{\partial^2}{\partial\mu_i^2}g$. There exists a linear transformation on $\mathbf{X}$ that converts the information matrix to the identity. Thus, the $L_{i,r}^{-1}$ factor is constant in the parameter.

THEOREM 2 (Admissibility). *For $p=1$ and $p=2$ the Jeffreys prior is admissible: there is no other prior $g$ such that*

$$\Delta g \le 0 \qquad \text{for all } \mu,$$
$$\Delta g < 0 \qquad \text{for some } \mu.$$

*For $p \ge 3$, the Jeffreys prior is inadmissible: there exists a prior $g$ such that $\Delta g < 0$ for all $\mu$; however, there exists no prior $g$ that dominates Jeffreys' uniformly by a positive amount, that is, there exists no $g$ such that for some $c > 0$,*

$$\frac{\Delta g}{g} \le -c < 0 \qquad \text{for all } \mu.$$

COROLLARY 2 (Minimaxity). *In one and two dimensions, the admissibility of the Jeffreys prior supports its minimaxity by the constant asymptotic risk. For location models of higher dimensionality, Jeffreys' is also minimax because it cannot be dominated uniformly.*

PROOF OF THEOREM 2. The case $p=1$ is a standard convex functions result. The case $p=2$ is simply Liouville's theorem (see [15]).

*Case* $p \ge 3$: If we assume that $\Delta g < 0$ everywhere, we can find a prior $g$ that makes the Jeffreys prior inadmissible. An example of this kind is

$$g(\boldsymbol{\mu}) = \left(1 + \sum_{i=1}^{p} \mu_i^2\right)^{\alpha},$$

which, for $0 > \alpha > 1 - \frac{p}{2}$, satisfies the condition $\Delta g < 0$ for any $\boldsymbol{\mu}$.

Following Aslan [2], we consider as our domain in $\mathcal{R}^p$ a solid sphere $D$ with radius $r$ and its surface $S:\{R=r\}$. We also consider the one-to-one mapping to "cylindrical" coordinates

$$\boldsymbol{\mu} \leftrightarrow (R, \mathbf{s}),$$

with $|\boldsymbol{\mu}| \le R^2$, $\mathbf{s}$ being of dimension $p-1$ and the Jacobian of the mapping being

$$\frac{\partial(\boldsymbol{\mu})}{\partial(R, \mathbf{s})} = R^{p-1}.$$

Gauss' divergence theorem applies and

$$\int \cdots \int_S r^{p-1} \frac{\partial g}{\partial r}\, d\mathbf{s} = \int \int \cdots \int_D R^{p-1} \Delta g \, dR \, d\mathbf{s},$$

or simply

$$r^{p-1} \bar{g}' = \int \int \cdots \int_D R^{p-1} \Delta g \, dR \, d\mathbf{s},$$

where $\bar{g}$ is the function obtained by averaging $g$ over $S$, $\bar{g}(r) = \int \cdots \int_S g(R, \mathbf{s})\, d\mathbf{s} \geq 0$ and $\bar{g}'$ is its derivative. By using the Leibniz rule for differentiation and the hypothesized inequality $\Delta g \leq -cg$ for a positive constant $c \in \mathcal{R}$, we obtain

$$\frac{d}{dr}(r^{p-1}\bar{g}') \leq -cr^{p-1}\bar{g}.$$

Making the change of variable $u = r^{2-p}$ with $du = (2-p)r^{1-p}\, dr$ and absorbing all the constants into $c$, we obtain the new differential inequality

$$\bar{g}'' \leq -cu^{\lambda}\bar{g} \qquad \forall\, u \geq 0, \tag{3.1}$$

where $-6 \leq \lambda = \frac{3(p-1)}{2-p} \leq -3$.

Assume first that $\bar{g}'(u) < 0$ for all $u$. This implies that $\bar{g}$ is strictly decreasing. Using the Taylor series expansion around $u_0 \leq u$ and (3.1), we obtain

$$\bar{g}(u) \leq \bar{g}(u_0) + (u - u_0)\bar{g}'(u_0) - \frac{cu^{\lambda}}{2}(u - u_0)^2\,\bar{g}(u^*),$$

where $u_0 \leq u^* \leq u$. Since $\bar{g}'(u_0) < 0$, $\bar{g}(u) \to -\infty$ as $u \to +\infty$, and this holds for any $u_0 \geq 0$. Therefore, $\bar{g}'$ must be nonnegative for all $u$.

Similarly, using a Taylor series expansion around $\frac{1}{2}u$ for $u \in [0, t]$, $t$ small, and by the strict concavity of $\bar{g}$, together with (3.1) and the increasing monotonicity of $\bar{g}$ in a small neighborhood of $u$, we obtain

$$\bar{g}(\tfrac{1}{2}u)[1 + \tfrac{1}{8}cu^{2+\lambda}] \leq \bar{g}(u).$$

Since $\lambda + 2 < -1$, we have $u^{2+\lambda} \to \infty$ as $u \to 0$. Thus, for small $u$, $\bar{g}(\frac{1}{2}u) \leq \frac{1}{2}\bar{g}(u)$, which contradicts the strict concavity of $\bar{g}$. The Jeffreys prior is inadmissible, but remains minimax for the case where the dimensionality of the model is at least three. □

In general, for invariant problems where the likelihood term in the asymptotic risk is constant, finding minimax solutions reduces to finding least favorable priors for which the prior term is minimax. If a reparametrization of the problem exists in which the Fisher information matrix is the identity, then the same argument used in the location case shows that the Jeffreys prior is asymptotically minimax.

## APPENDIX: ASYMPTOTICS OF RISK

We present here the main assumptions and ideas used in obtaining the result of Theorem 1. For a more elaborate and computationally involved presentation of the proof, see [2].

The asymptotics of the risk involve both Taylor series and Edgeworth approximations which require appropriate regularity conditions. The work of Bhattacharya and Ghosh [3] gives a rigorous account of the theory of Edgeworth series for general statistics. The Taylor series approximations are polynomials in $(\hat{\boldsymbol{\theta}} - \boldsymbol{\theta})$ with remainder terms which require special attention in order to integrate successfully. We also need expectations to evaluate risks accurately to $O(n^{-2})$ terms.

The locally asymptotic normality of the standardized maximum likelihood estimator $\sqrt{n}(\hat{\boldsymbol{\theta}} - \boldsymbol{\theta})$, as well as truncated expectations in the sense given by Hartigan [8] is used extensively here. The following assumptions, similar to Hartigan's, are required for the validity of our expansions:

Assumptions. (A1) The prior density $h$ is smooth in the sense that it is twice continuously differentiable in a neighborhood of $\boldsymbol{\theta}$ and is positive.

(A2) We assume that the second derivatives $\frac{\partial^2}{\partial \theta_i^2} l(\boldsymbol{\theta})$ are of order $n$, where $n$ is the number of observations. We also assume that all $l$'s and $L$'s are, in general, of order $n$ with the exception of $l_i$'s, which are random variables of order $\sqrt{n}$ with zero expectations. These assumptions are usually satisfied in practice.

(A3) $l(\boldsymbol{\theta}) = \log f(\mathbf{x}|\boldsymbol{\theta})$ is five times continuously differentiable with respect to $\tilde{\boldsymbol{\theta}}$ in a small neighborhood of the true parameter $\boldsymbol{\theta}$ for each observation $\mathbf{x}$.

(A4) All moments exist for the first four log-likelihood derivatives, and for the maximum squared fifth derivatives, in a neighborhood of $\boldsymbol{\theta}$; in other words, for each $\tilde{\boldsymbol{\theta}} \in \Theta$ and for some $\varepsilon > 0$,

$$\mathbf{P}_{\tilde{\boldsymbol{\theta}}} |l_{i_1,\ldots,i_4}(\tilde{\boldsymbol{\theta}})| < \infty,$$

$$\mathbf{P}_{\tilde{\boldsymbol{\theta}}} \Bigl( \sup_{|\boldsymbol{\theta} - \tilde{\boldsymbol{\theta}}| \le \varepsilon} |l_{i_1,\ldots,i_5}(\boldsymbol{\theta})| \Bigr)^2 < \infty,$$

where the set of indices $(i_1, \ldots, i_r) \subseteq \{1, \ldots, p\}^r$, with $r = 4$ and $r = 5$, respectively.

(A5) The integral $\int f(\mathbf{x}|\boldsymbol{\theta})\, d\mathbf{x}$ can be differentiated four times with respect to $\boldsymbol{\theta}$ under the integral sign. The usual likelihood identities, obtained by differentiating $\int f_{\boldsymbol{\theta}}$, are valid for any indices $i$, $j$, $k$ and $l$:

$$\begin{aligned}
0 &= L_i, \\
0 &= L_{ij} + L_{i,j}, \\
0 &= L_{ijk} + L_{ij,k} + L_{ik,j} + L_{jk,i} + L_{i,j,k}, \\
0 &= L_{ijkl} + L_{ijk,l} + L_{ijl,k} + L_{ikl,j} \\
&\quad + L_{jkl,i} + L_{ij,kl} + L_{ik,jl} + L_{il,jk} \\
&\quad + L_{ij,k,l} + L_{ik,j,l} + L_{il,j,k} \\
&\quad + L_{jk,i,l} + L_{jl,i,k} + L_{kl,i,j} + L_{i,j,k,l}.
\end{aligned}$$

(A6) The Fisher information matrix $(-L_{ij})_{i,j=1,\dots,p} = (L_{i,j})_{i,j=1,\dots,p}$, or $L_{i,j}$ for short, is nonsingular and positive definite for $|\boldsymbol{\theta} - \tilde{\boldsymbol{\theta}}| < \varepsilon$.

(A7) For each $\varepsilon > 0$, $\mathbf{P}\{|\boldsymbol{\theta} - \tilde{\boldsymbol{\theta}}| > \varepsilon\} = o(n^{-2})$.

PROOF OF THEOREM 1. Through straightforward calculations the risk expression becomes a difference between two Kullback–Leibler losses. When $\boldsymbol{\theta}$ is true, we have

$$R(\boldsymbol{\theta}, f_h) = D(f(\mathbf{x_{n+1}}|\boldsymbol{\theta})||f(\mathbf{x_{n+1}})) - D(f(\mathbf{x_n}|\boldsymbol{\theta})||f(\mathbf{x_n})), \tag{A.1}$$

where $f(\mathbf{x})$ stands for the marginal density of $\mathbf{x}$.

To arrive at the risk asymptotics, we begin by computing an asymptotic expression for the Kullback–Leibler loss,

$$D(f(\mathbf{x}|\boldsymbol{\theta})||f(\mathbf{x})) = \int \{\log f(\mathbf{x}|\boldsymbol{\theta}) - \log f(\mathbf{x})\} f(\mathbf{x}|\boldsymbol{\theta})\, d\boldsymbol{\theta}.$$

The following lemma gives the asymptotic behavior of $f(\mathbf{x})$. The result and its proof can be found in [2].

LEMMA 2. *Under the previous regularity conditions, the marginal density $f(\mathbf{x})$ has the following asymptotic expression with terms of order $O_P(n^{-1})$, ignoring smaller $O_P(n^{-2})$ terms:*

$$\begin{aligned} f(\mathbf{x}) = {} & (2\pi)^{p/2} \det(-\hat{l}_{ij})^{-1/2} f(\mathbf{x}|\hat{\boldsymbol{\theta}}) h(\hat{\boldsymbol{\theta}}) \\ & \times \{1 - \tfrac{1}{2} L_{ij}^{-1} \hat{h}_{ij} - \tfrac{1}{2} L_{ii'}^{-1} \hat{h}_i \hat{h}_{i'} \\ & + \tfrac{1}{8} L_{ij}^{-1} L_{kl}^{-1} \hat{l}_{ijkl} - \tfrac{1}{12} L_{ii'}^{-1} L_{jj'}^{-1} L_{kk'}^{-1} \hat{l}_{ijk} \hat{l}_{i'j'k'} \\ & - \tfrac{1}{8} L_{ij}^{-1} L_{i'j'}^{-1} L_{kk'}^{-1} \hat{l}_{ijk} \hat{l}_{i'j'k'} + \tfrac{1}{2} L_{ij}^{-1} L_{ki'}^{-1} \hat{l}_{ijk} \hat{h}_{i'} + O_P(n^{-2})\}. \end{aligned}$$

Using Lemma 2 and the "Expectation lemma" in [8], the Kullback–Leibler loss expression becomes

$$\begin{aligned} D(f(\mathbf{x}|\boldsymbol{\theta})||f(\mathbf{x})) = {} & -\frac{p}{2} \log 2\pi + \int (\log f(\mathbf{x}|\boldsymbol{\theta}) - \log f(\mathbf{x}|\hat{\boldsymbol{\theta}})) f(\mathbf{x}|\boldsymbol{\theta})\, d\mathbf{x} \\ & - \int \log \frac{h(\hat{\boldsymbol{\theta}})}{\det(-\hat{l}_{ij})^{1/2}} f(\mathbf{x}|\boldsymbol{\theta})\, d\mathbf{x} \\ & + \Big\{ \frac{1}{2} L_{ij}^{-1} h_{ij} + \frac{1}{2} L_{ii'}^{-1} h_i h_{i'} \\ & - \frac{1}{8} L_{ij}^{-1} L_{kl}^{-1} L_{ijkl} + \frac{1}{12} L_{ii'}^{-1} L_{jj'}^{-1} L_{kk'}^{-1} L_{ijk} L_{i'j'k'} \\ & + \frac{1}{8} L_{ij}^{-1} L_{i'j'}^{-1} L_{kk'}^{-1} L_{ijk} L_{i'j'k'} \\ & - \frac{1}{2} L_{ij}^{-1} L_{ki'}^{-1} L_{ijk} h_{i'} + O_P(n^{-2}) \Big\}. \end{aligned}$$

Following Aslan [2], the first two integrals in the Kullback–Leibler expression are further expanded into the following asymptotic expressions:

$$\int (\log f(\mathbf{x}|\boldsymbol{\theta}) - \log f(\mathbf{x}|\hat{\boldsymbol{\theta}}))f(\mathbf{x}|\boldsymbol{\theta})\,d\mathbf{x}$$

$$= -\frac{p}{2} + L_{i,r}^{-1}L_{j,s}^{-1}\Big(-\frac{1}{2}L_{i,j,rs} - \frac{1}{2}L_{ij,rs} - \frac{1}{2}L_{irj,s} - \frac{1}{8}L_{irjs}\Big)$$

$$+ L_{i,r}^{-1}L_{j,s}^{-1}L_{k,t}^{-1}\Big(-\frac{1}{2}L_{i,rj}L_{k,st} - \frac{1}{2}L_{i,st}L_{k,rj}$$

$$- \frac{1}{6}L_{ijk}L_{r,s,t} - \frac{1}{2}L_{irj}L_{k,st}$$

$$- L_{ijk}L_{r,st} - \frac{1}{6}L_{irj}L_{skt} - \frac{1}{4}L_{ijk}L_{rst}\Big) + O_P(n^{-2})$$

and

$$-\int \log \frac{h(\hat{\boldsymbol{\theta}})}{\det(-\hat{l}_{ij})^{1/2}} f(\mathbf{x}|\boldsymbol{\theta})d\mathbf{x}$$

$$= \frac{p}{4} + \frac{1}{2}\log(|L_{i,j}|) - \log h$$

$$+ L_{i,r}^{-1}L_{j,s}^{-1}\Big(-\frac{1}{2}L_{irj,s} - \frac{1}{4}L_{ij,rs} - \frac{1}{4}L_{irjs}\Big)$$

$$+ L_{i,r}^{-1}L_{j,s}^{-1}L_{k,t}^{-1}\Big(-\frac{1}{2}L_{irj}L_{k,st} - \frac{1}{4}L_{irj}L_{skt} - \frac{1}{2}L_{ijk}L_{r,st} - \frac{1}{4}L_{ijk}L_{rst}\Big)$$

$$+ L_{i,r}^{-1}L_{j,s}^{-1}\Big(-L_{rj,s}h_i - \frac{1}{2}L_{rjs}h_i\Big) + \frac{1}{2}L_{i,r}^{-1}h_{ir} + O_P(n^{-2}).$$

By simply substituting these expressions into the Kullback–Leibler loss formula from above, we obtain the following asymptotic approximation of the loss:

LEMMA 3. *Under the previous regularity conditions, the asymptotic expression for the Kullback–Leibler loss function $D(f(\mathbf{x}|\boldsymbol{\theta})||f(\mathbf{x}))$ with terms of order $O_P(n^{-1})$ and ignoring smaller terms of order $O_P(n^{-2})$ is as follows:*

$$D(f(\mathbf{x}|\boldsymbol{\theta})||f(\mathbf{x}))$$

$$= -\frac{p}{2}\log 2\pi - \frac{p}{4} + \frac{1}{2}\log(|L_{i,j}|) - \log h$$

$$+ L_{i,r}^{-1}L_{j,s}^{-1}\Big(-\frac{1}{2}L_{ij,r,s} - \frac{3}{4}L_{ij,rs} - L_{irj,s} - \frac{1}{2}L_{irjs}\Big)$$

$$+ L_{i,r}^{-1}L_{j,s}^{-1}L_{k,t}^{-1}\Big(-\frac{1}{2}L_{i,rj}L_{k,st} - \frac{1}{2}L_{i,jk}L_{t,rs} - \frac{1}{6}L_{ijk}L_{r,s,t}$$

$$- L_{irj}L_{k,st} - \frac{3}{2}L_{ijk}L_{r,st} - \frac{1}{2}L_{irj}L_{skt} - \frac{7}{12}L_{ijk}L_{rst}\Big)$$

$$+ L_{i,r}^{-1}L_{j,s}^{-1}(-L_{rj,s} - L_{rjs})h_i + L_{i,r}^{-1}\left(-h_{ir} - \frac{1}{2}h_i h_r\right)$$

$$+ O_P(n^{-2}).$$

We now arrive at the asymptotic risk approximation in Theorem 1 by simply substituting in (A.1) the two asymptotic Kullback–Leibler loss expressions, written more concisely as

$$D_{n+1} = -\frac{p}{2}\log 2\pi - \frac{p}{4} + \frac{1}{2}\log\{(n+1)^p|L_{i,j}|\} - \log h(\boldsymbol{\theta}) - \frac{G(\boldsymbol{\theta})}{n+1} + O_P((n+1)^{-2})$$

and

$$D_n = -\frac{p}{2}\log 2\pi - \frac{p}{4} + \frac{1}{2}\log\{n^p|L_{i,j}|\} - \log h(\boldsymbol{\theta}) - \frac{G(\boldsymbol{\theta})}{n} + O_P(n^{-2}),$$

where $D_{n+1} = D(f(\mathbf{x_{n+1}}|\boldsymbol{\theta})||f(\mathbf{x_{n+1}}))$ and $D_n = D(f(\mathbf{x_n}|\boldsymbol{\theta})||f(\mathbf{x_n}))$. Thus, the difference $D_{n+1} - D_n$ will be of the form

$$\frac{p}{2n} - \frac{p}{4n^2} + \frac{G(\boldsymbol{\theta})}{n^2} + O_P(n^{-3}),$$

where $G(\boldsymbol{\theta})$ is the $n\{\cdots\}$ term in the asymptotic risk expression of Theorem 1. $\square$

**Acknowledgments.** This work represents a portion of my doctoral dissertation, written at Yale University under the supervision of John Hartigan. I am very grateful to him for his generous insights and suggestions. I would also like to thank the referee and the Associate Editor for their useful comments that significantly improved this presentation.

Clinical Epidemiology Research Center
VA Connecticut Healthcare System
950 Campbell Avenue (Mailcode 151B)
West Haven, Connecticut 06516
USA
E-mail: mihaela.aslan@yale.edu